\documentclass[a4paper,10pt]{article}
\usepackage[dvips]{graphicx}
\usepackage[T1]{fontenc}
\usepackage{times}
\usepackage{mltex}
\usepackage{amsmath,amsthm,amsfonts,amssymb}
\usepackage{pstcol,pst-fill,pst-grad,pst-3d,pst-plot}
\usepackage{pstricks}
\usepackage{subfigure}
\usepackage{cases}
\usepackage{shorttoc}
\usepackage{fancyhdr,calc}
\usepackage{multicol}
\usepackage[all]{xy}
\usepackage{fancyhdr}

\newcommand{\R}{\mathbb{R}}
\newcommand{\C}{\mathbb{C}}

\begin{document}

\title{Extracting squares from any quadrilateral}

\author{Pierre Godard}

\date{}

\maketitle

\begin{abstract}
A procedure that generates parallelograms from any quadrilateral is presented. If the original quadrilateral is itself a parallelogram, then the procedure gives squares. Hence, when applied two times, this procedure generates squares from any quadrilateral. The proof needs only undergraduate level.
\end{abstract}

\vspace{.8cm}

\section{Notations}
We identify the plane $\R^2$ with the complex plane $\C$; the imaginary unit is denoted by $i$.
The tools this article needs are very basic and may be found in any textbooks on geometry, \emph{e.g.} the first 9 chapters of \cite{bCoxeter69}.

Let $r_{\alpha}$ be the rotation about the point $O$ by angle $\alpha$, and let $t_{\overrightarrow{Oa}}$ be the translation by the vector $\overrightarrow{Oa}$.
Then the rotation about the point $a$ by angle $\alpha$ is
\begin{equation*}
r_{a,\alpha}=t_{\overrightarrow{Oa}}\circ r_{\alpha}\circ(t_{\overrightarrow{Oa}})^{-1}:z\mapsto e^{i\alpha}(z-a)+a.
\end{equation*}
Composition of operations, denoted here by $\circ$, will no more be explicitely written.

We use the notation $P=[a_1,a_2,a_3,a_4]$ when the four vertices of the quadrilateral $P$ are $a_1$, $a_2$, $a_3$ and $a_4$, and the edges are $a_1a_2$, $\dots$, $a_4a_1$.
Note that the symbols can be circularly permuted or the order can be reversed;
for example, we have $[a_1,a_2,a_3,a_4]=[a_2,a_3,a_4,a_1]=[a_2,a_1,a_4,a_3]$ but $[a_1,a_2,a_3,a_4]\neq[a_1,a_3,a_2,a_4]$.

\section{Exctracting parallelograms from any quadrilateral}
\newtheorem{th_parallelogram}{theorem}
\begin{th_parallelogram}{Let $P$ be a quadrilateral and denote its four vertices by $a_1$, $a_2$, $a_3$ and $a_4$; let $\alpha_i$ be the half of the angle at the vertex $a_i$. Moreover, let $b_{ijkl}$ be the fixed point of the rotation $r_{a_i,\alpha_i}r_{a_j,\alpha_j}r_{a_k,\alpha_k}r_{a_l,\alpha_l}$, with $\{i,j,k,l\}=\{1,2,3,4\}$.
Then the following sets of points define parallelograms:
}\label{th_parallelogram1}

\begin{align*}
P^e:&=[b_{1234},b_{1243},b_{2143},b_{2134}] \\
P^{(13)(24)}:&=[b_{3412},b_{3421},b_{4321},b_{4312}] \\
P^{(23)}:&=[b_{1324},b_{1342},b_{3142},b_{3124}] \\
P^{(1342)}:&=[b_{2413},b_{2431},b_{4231},b_{4213}] \\
P^{(24)}:&=[b_{1432},b_{1423},b_{4123},b_{4132}] \\
P^{(13)}:&=[b_{3214},b_{3241},b_{2341},b_{2314}].
\end{align*}
Moreover, $P^e$ is congruent to $P^{(13)(24)}$, $P^{(23)}$ is congruent to $P^{(1342)}$, $P^{(24)}$ is congruent to $P^{(13)}$.
\end{th_parallelogram}

First, we note that $b_{ijkl}$ is well-defined since $\alpha_i+\alpha_j+\alpha_k+\alpha_l=\pi$ when $i$, $j$, $k$ and $l$ are all different.
The proof of the theorem proceeds in three steps:
\newtheorem{lemma_parallelogram_1}{lemma}
\begin{lemma_parallelogram_1}{Under the hypotheses of the theorem \ref{th_parallelogram1},
$P^e$ is a parallelogram.
}\label{th_lemma1}
\end{lemma_parallelogram_1}

\begin{proof}
For convenience, we write $r_1$ for $r_{a_1,\alpha_1}$, and similarly for the other points of $P$.
$r_4r_3r_2r_1$ is the rotation about $b_{4321}$ by the angle $\pi$.
Hence it is an involution.
Similarly, $r_1r_2r_3r_4$ is an involution.
Thus, we have
\begin{align*}
(r_1r_2r_3r_4r_1r_2r_4r_3)(r_2r_1r_4r_3r_2r_1r_3r_4)&=r_1r_2r_3r_4r_1r_2(r_4r_3r_2r_1r_4r_3r_2r_1)r_3r_4 \\
&=r_1r_2r_3r_4r_1r_2r_3r_4 \\
&=id
\end{align*}
where $id:z\mapsto z$.
On the other hand
\begin{align*}
(r_1r_2r_3r_4r_1r_2r_4r_3)(r_2r_1r_4r_3r_2r_1r_3r_4)&=(r_{b_{1234},\pi}r_{b_{1243},\pi})(r_{b_{2143},\pi}r_{b_{2134},\pi}) \\
&=t_{2\overrightarrow{b_{1243}b_{1234}}}t_{2\overrightarrow{b_{2134}b_{2143}}} \\
&=t_{2(\overrightarrow{b_{1243}b_{1234}}+\overrightarrow{b_{2134}b_{2143}})}.
\end{align*}
Hence $\overrightarrow{b_{1243}b_{1234}}+\overrightarrow{b_{2134}b_{2143}}=0$,
or $[b_{1243},b_{1234},b_{2134},b_{2143}]$ is a parallelogram.
\end{proof}

\begin{lemma_parallelogram_1}{Under the hypotheses of the theorem \ref{th_parallelogram1}, if $P^e$ is a parallelogram, then
$P^{(13)(24)}$, $P^{(23)}$, $P^{(1342)}$, $P^{(24)}$ and $P^{(13)}$ are all parallelograms.
}\label{th_lemma2}
\end{lemma_parallelogram_1}

\begin{proof}
There is an action of the symmetric group on four symbols $S_4$ to the four points of $P$.
This induces an action of $S_4$ on the $b_{ijkl}$'s: for $\sigma$ in $S_4$, we have
\begin{equation*}
\sigma:b_{ijkl}\mapsto b_{\sigma^{-1}(i)\sigma^{-1}(j)\sigma^{-1}(k)\sigma^{-1}(l)}.
\end{equation*}
When $S_4$ acts on the identity $\overrightarrow{b_{1243}b_{1234}}+\overrightarrow{b_{2134}b_{2143}}=0$ deduced in lemma \ref{th_lemma1}, other parallelograms than $P^e$ are evidenced.
Obviously, the group $C_2\times C_2$ generated by the permutations $(12)$ and $(34)$ leaves the vertex set of $P^e$ invariant.
Now we choose as representatives of the left cosets of $S_4/(C_2\times C_2)$ the elements $e$, $(13)(24)$, $(23)$, $(1342)$, $(24)$ and $(13)$.
The action of these elements on the set $\{b_{1234},b_{1243},b_{2143},b_{2134}\}$ gives respectively the vertex set of $P^e$, $P^{(13)(24)}$, $P^{(23)}$, $P^{(1342)}$, $P^{(24)}$ and $P^{(13)}$.
\end{proof}
 
\begin{lemma_parallelogram_1}{Under the hypotheses of the theorem \ref{th_parallelogram1},
$P^e$ is congruent to $P^{(13)(24)}$, $P^{(23)}$ is congruent to $P^{(1342)}$, $P^{(24)}$ is congruent to $P^{(13)}$.
}\label{th_lemma3}
\end{lemma_parallelogram_1}

\begin{proof}
The proof is straightforward: we first compute that
\begin{align*}
&b_{ijkl}=\frac{1}{2}\Big(a_i(1-e^{i\alpha_i})+a_je^{i\alpha_i}(1-e^{i\alpha_j}) \\
&\qquad+a_ke^{i(\alpha_i+\alpha_j)}(1-e^{i\alpha_k})+a_le^{i(\alpha_i+\alpha_j+\alpha_k)}(1-e^{i\alpha_l})\Big)
\end{align*}
when $\{i,j,k,l\}=\{1,2,3,4\}$.
Then, we easily check that
\begin{align*}
b_{1234}-b_{1243}&=\frac{1}{2}e^{i(\alpha_1+\alpha_2)}(1-e^{i\alpha_3})(1-e^{i\alpha_4})(a_3-a_4) \\
&=e^{i(\alpha_1+\alpha_2)}(b_{3421}-b_{4321}) \\
b_{1234}-b_{2134}&=\frac{1}{2}(1-e^{i\alpha_1})(1-e^{i\alpha_2})(a_1-a_2) \\
&=e^{i(\alpha_1+\alpha_2)}(b_{3421}-b_{3412})
\end{align*}
so that $P^e$ and $P^{(13)(24)}$ are congruent.
The action of $(23)$ (respectively $(24)$) on these equalities implies that $P^{(23)}$ is congruent to $P^{(1342)}$ (respectively that $P^{(24)}$ is congruent to $P^{(13)}$).
\end{proof}

\begin{figure}
\centering
\includegraphics[width=\textwidth]{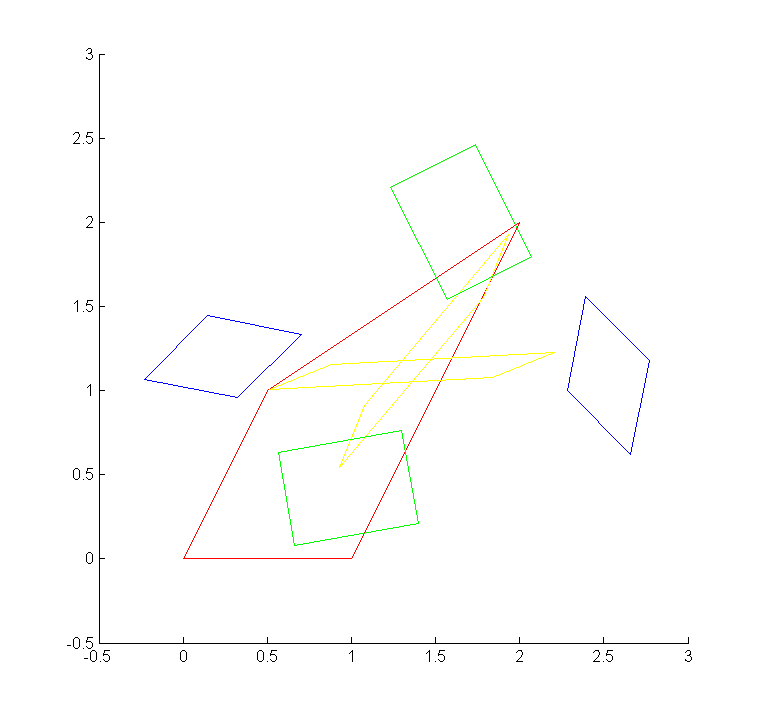}
\caption{Example illustrating the theorem \ref{th_parallelogram1}. The coordinates of $a_1$, $a_2$, $a_3$ and $a_4$ are respectively
$(0,0)$, $(1,0)$, $(2,2)$ and $(1/2,1)$. The quadrilateral $P$ is red, the parallelograms $P^e$ and $P^{(13)(24)}$ are blue, the parallelograms $P^{(23)}$ and $P^{(1342)}$ are yellow and the parallelograms $P^{(24)}$ and $P^{(13)}$ are green.}
\label{fig_parallelogram}
\end{figure}

The figure \ref{fig_parallelogram} presents an example of the quadrilateral $P$ and the six unveiled parallelograms.

\section{Exctracting squares from any parallelogram}
With the notation of the theorem \ref{th_parallelogram1}, we note that if $P$ is a parallelogram then
\begin{equation}
     \begin{cases}
	\alpha_1=\alpha_3 \\
	\alpha_2=\alpha_4 \\
	a_1-a_2=-(a_3-a_4) \\
	e^{i(\alpha_1+\alpha_2)}=i.
     \end{cases}
\label{cond_parallelogram}
\end{equation}
Hence, from the proof of the lemma \ref{th_lemma3}, $b_{1234}-b_{1243}=-i(b_{1234}-b_{2134})$,
or the quadrilateral $P^e$ has two adjacent edges of the same size that form a right angle.
Since we know from the theorem \ref{th_parallelogram1} that $P^e$ is a parallelogram, it is a square.
The set of conditions \ref{cond_parallelogram} is invariant with respect to the permutations $(13)(24)$,
and an other set of valid conditions is obtained with the actions of $(24)$ or $(13)$.
On the contrary, the actions of $(23)$ or $(1342)$ are not allowed since $a_1-a_3\neq-(a_2-a_4)$, the segments $a_2a_4$ and $a_1a_3$ being the diagonals of a parallelogram.
The proof of the lemma \ref{th_lemma3} also shows that if $P$ is a parallelogram, then
\begin{equation*}
b_{1234}-b_{1243}=i(b_{3421}-b_{4321}),
\end{equation*}
so that $P^e$ and $P^{(13)(24)}$ are equal up to a translation.
The action of $(24)$ on this last equation implies that $P^{(24)}$ is equal to $P^{(13)}$, once again up to a translation.
We thus have proved the following theorem:

\begin{th_parallelogram}{Let $P$ be a parallelogram and denote its four vertices by $a_1$, $a_2$, $a_3$ and $a_4$; let $\alpha$ be the half of the angle at the vertex $a_1$. Moreover, let $b_{ijkl}$ be the fixed point of the rotation $r_{a_i,\alpha_i}r_{a_j,\alpha_j}r_{a_k,\alpha_k}r_{a_l,\alpha_l}$, with $i$, $j$, $k$ and $l$ all different.
Then the following sets of points define squares:
}\label{th_parallelogram2}

\begin{align*}
P^e:&=[b_{1234},b_{1243},b_{2143},b_{2134}] \\
P^{(13)(24)}:&=[b_{3412},b_{3421},b_{4321},b_{4312}] \\
P^{(24)}:&=[b_{1432},b_{1423},b_{4123},b_{4132}] \\
P^{(13)}:&=[b_{3214},b_{3241},b_{2341},b_{2314}].
\end{align*}
Moreover, $P^e$ can be transformed into $P^{(13)(24)}$ by a translation, and $P^{(24)}$ can be transformed into $P^{(13)}$ by a translation.
Finally, the following sets of points define parallelograms:
\begin{align*}
P^{(23)}:&=[b_{1324},b_{1342},b_{3142},b_{3124}] \\
P^{(1342)}:&=[b_{2413},b_{2431},b_{4231},b_{4213}].
\end{align*}
\end{th_parallelogram}

\begin{figure}
\centering
\includegraphics[width=\textwidth]{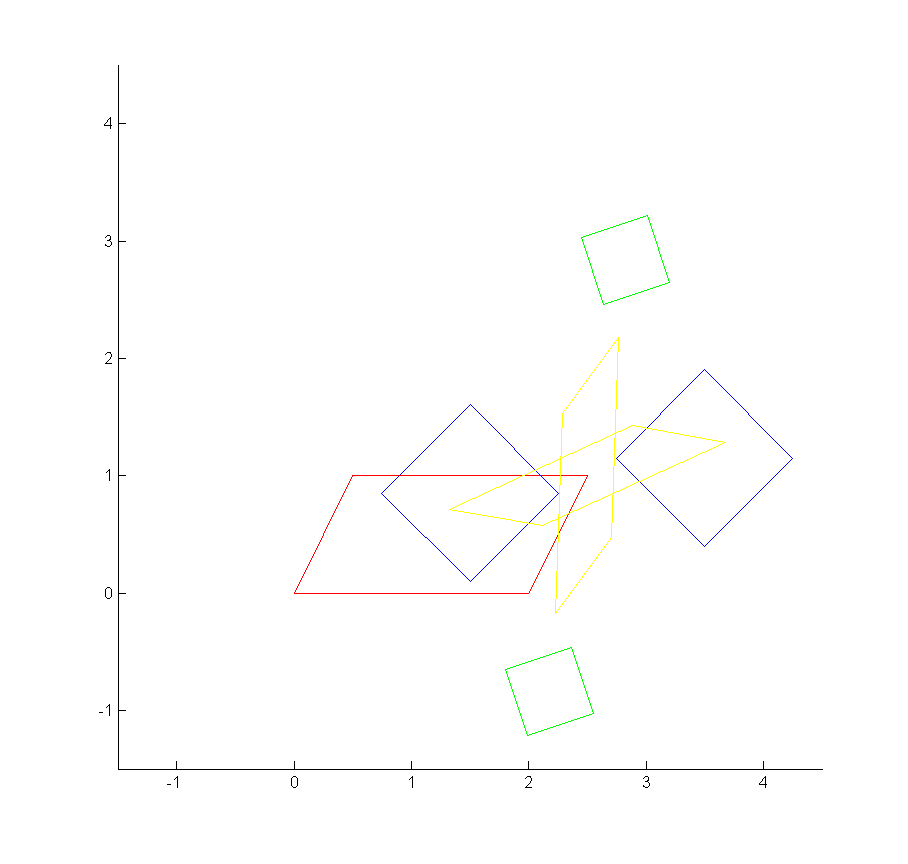}
\caption{Example illustrating the theorem \ref{th_parallelogram2}. The coordinates of $a_1$, $a_2$, $a_3$ and $a_4$ are respectively
$(0,0)$, $(2,0)$, $(5/2,1)$ and $(1/2,1)$. The parallelogram $P$ is red, the squares $P^e$ and $P^{(13)(24)}$ are blue, the parallelograms $P^{(23)}$ and $P^{(1342)}$ are yellow and the squares $P^{(24)}$ and $P^{(13)}$ are green.}
\label{fig_square}
\end{figure}

That $P^{(23)}$ and $P^{(1342)}$ are parallelograms is of course a direct consequence of the theorem \ref{th_parallelogram1}.
The figure \ref{fig_square} presents an example of the parallelogram $P$ and the four unveiled squares;
the two parallelograms $P^{(23)}$ and $P^{(1342)}$ are also represented.

\begin{figure}
\centering
\includegraphics[width=\textwidth]{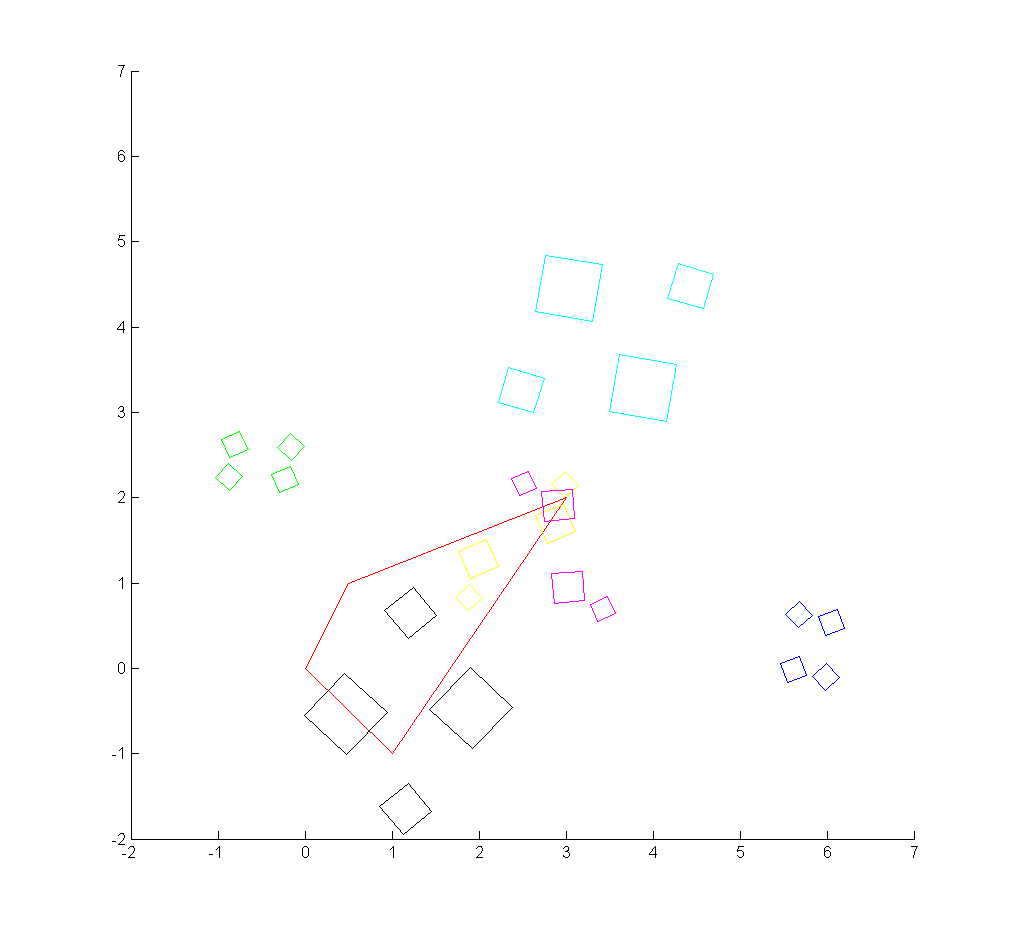}
\caption{Example illustrating a combination of the theorem \ref{th_parallelogram1} and the theorem \ref{th_parallelogram2}. The coordinates of $a_1$, $a_2$, $a_3$ and $a_4$ are respectively
$(0,0)$, $(1,-1)$, $(3,2)$ and $(1/2,1)$. The parallelogram $P$ is red; $P$ gives rise to six (not represented) parallelograms, $P^e$, $P^{(13)(24)}$, $P^{(23)}$, $P^{(1342)}$, $P^{(24)}$ and $P^{(13)}$, each of which gives rise to four squares. The squares coming from $P^e$, $P^{(13)(24)}$, $P^{(23)}$, $P^{(1342)}$, $P^{(24)}$ and $P^{(13)}$ are respectively colored in blue, green, yellow, magenta, black and cyan.}
\label{fig_24square}
\end{figure}

\section{Conclusion and outlook}
Starting from any quadrilateral $P$, the theorem \ref{th_parallelogram1} gives six parallelograms.
For each of these parallelograms, the theorem \ref{th_parallelogram2} gives four squares.
As illustrated in the figure \ref{fig_24square}, which presents an example of a quadrilateral $P$ and the unveiled squares,
these twenty-four squares can all be different.

The two sets of equations presented in the lemma \ref{th_lemma1} can be generalized to any polygons with an even number of edges.
For example, if $a_1$, $a_2$, $a_3$, $a_4$, $a_5$ and $a_6$ denote the six points of a hexagon and $\alpha_i$ is one fourth of the angle at the vertex $a_i$, then
\begin{equation*}
r_1r_2r_3r_4r_5r_6r_1r_2r_3r_5r_6r_4r_2r_3r_1r_5r_6r_4r_2r_3r_1r_6r_4r_5r_3r_1r_2r_6r_4r_5r_3r_1r_2r_4r_5r_6=id,
\end{equation*}
with $r_i=r_{a_i,\alpha_i}$, which means that
\begin{equation*}
\overrightarrow{b_{123564}b_{123456}}+\overrightarrow{b_{231645}b_{231564}}+\overrightarrow{b_{312456}b_{312645}}=0.
\end{equation*}
Hence, the hexagon $[b_{123564},b_{123456},b_{231645},b_{231564},b_{312456},b_{312645}]$ has not twelve degrees of freedom, but only ten.
We could for example inquire about a procedure which, when applied a sufficient number of times, eventually leads to a regular hexagon.

\vspace{1cm}

\begin{flushright}
\emph{The author gratefully aknowledges fruitful dscussions with Oleg Ogievetsky.}
\end{flushright}

\bibliographystyle{unsrt}

\bibliography{biblio}

\end{document}